\documentclass[draft]{elsart1p}

\usepackage{amsmath,amssymb,mathrsfs}
\usepackage{pb-diagram,lamsarrow,pb-lams}

\input xy
\xyoption{all}

\newcommand{\dmo}{\DeclareMathOperator*}
\dmo{\colim}{colim}

\newcommand{\ncb}[2]{\newcommand{#1}{\mathbb{#2}}}
\ncb{\CP}{CP} \ncb{\F}{F} \ncb{\bL}{L} \ncb{\N}{N} \ncb{\RP}{RP} \ncb{\Q}{Q} \ncb{\Z}{Z}

\newcommand{\ncc}[2]{\newcommand{#1}{\mathcal{#2}}}
\ncc{\cN}{N} \ncc{\cH}{H} \ncc{\cO}{O} \ncc{\cS}{S}

\newcommand{\ncf}[2]{\newcommand{#1}{\mathscr{#2}}}
\ncf{\fB}{B} \ncf{\fC}{C} \ncf{\fL}{L} \ncf{\fN}{N}

\newcommand{\nck}[2]{\newcommand{#1}{\mathfrak{#2}}}
\nck{\kp}{p}

\newcommand{\ncs}[2]{\newcommand{#1}{\boldsymbol{\mathsf{#2}}}}
\ncs{\sC}{C} \ncs{\sD}{D} \ncs{\sSD}{SD} \ncs{\sQ}{Q}

\newcommand{\nc}[2]{\newcommand{#1}{#2}}
\nc{\ABELIANGROUPS}{\mathrm{ABELIAN}\,\mathrm{GROUPS}}
\nc{\bdry}{\partial}
\nc{\COMMUTATIVERINGS}{\mathrm{COMMUTATIVE}\,\mathrm{RINGS}}
\nc{\DS}{\displaystyle} \nc{\en}{\enspace}
\nc{\FINITEGROUPS}{\mathrm{FINITE}\,\mathrm{GROUPS}}
\nc{\hookra}{\hookrightarrow} \nc{\inv}{^{-1}} \nc{\iso}{\cong}
\nc{\longra}{\longrightarrow}
\nc{\ST}{\;|\;} \nc{\vphi}{\varphi} \nc{\x}{\times}
\nc{\xo}{\otimes}

\newcommand{\ncm}[2]{\newcommand{#1}{\mathrm{#2}}}
\ncm{\Ad}{Ad} \ncm{\Arf}{Arf} \ncm{\aug}{aug} \ncm{\Aut}{Aut}
\ncm{\Cok}{Cok} \ncm{\Conj}{conj} \ncm{\diag}{diag} \ncm{\End}{End}
\ncm{\eval}{eval} \ncm{\hAut}{hAut} \ncm{\Hom}{Hom} \ncm{\id}{id}
\ncm{\incl}{incl} \ncm{\Ker}{Ker} \ncm{\Mor}{Mor} \ncm{\Ob}{Ob}
\ncm{\quotient}{quotient} \ncm{\Nil}{Nil} \ncm{\op}{op}
\ncm{\Split}{split} \ncm{\Tor}{Tor} \ncm{\UNil}{UNil}

\newcommand{\abs}[1]{\left| #1 \right|}
\newcommand{\gens}[1]{\left\langle #1 \right\rangle}
\newcommand{\ol}[1]{\overline{#1}}
\newcommand{\prn}[1]{\left( #1 \right)}
\newcommand{\set}[1]{\left\{\, #1 \,\right\}}
\newcommand{\SmMatrix}[1]{\left(\begin{smallmatrix} #1 \end{smallmatrix}\right)}
\newcommand{\wh}[1]{{\widehat{#1}}}
\newcommand{\wt}[1]{{\widetilde{#1}}}
\newcommand{\xra}[1]{\xrightarrow{#1}}

\newcommand{\infdec}{^{\gens{-\infty}}}

\begin{document}

\begin{frontmatter}
\title{Reduction of $\UNil$ for finite groups\\ with normal abelian Sylow 2-subgroup}
\author{Qayum Khan}
\address{Department of Mathematics, Vanderbilt University, Nashville, TN 37240
U.S.A.} \ead{qayum.khan@vanderbilt.edu}

\begin{abstract}
Let $F$ be a finite group with a Sylow 2-subgroup $S$ that is normal
and abelian. Using hyperelementary induction and cartesian squares,
we prove that Cappell's unitary nilpotent groups
$\UNil_*(\Z[F];\Z[F],\Z[F])$ have an induced isomorphism to the
quotient of $\UNil_*(\Z[S];\Z[S],\Z[S])$ by the action of the group
$F/S$. In particular, any finite group $F$ of odd order has the same
$\UNil$-groups as the trivial group. The broader scope is the study
of the $L$-theory of virtually cyclic groups, based on the
Farrell--Jones isomorphism conjecture. We obtain partial information
on these $\UNil$ when $S$ is a finite abelian 2-group and when $S$
is a special 2-group.
\end{abstract}

\end{frontmatter}

\section{Introduction}

Our main theorem reduces the computation of $\UNil$ for finite
groups with normal abelian Sylow 2-subgroup to the computation of
$\UNil$ for its Sylow 2-subgroup.  Throughout the paper, all
multiplicative groups are equipped with trivial orientation
character.

\begin{thm}\label{Thm_FiniteReduction}
Suppose $F$ is a finite group with a normal abelian Sylow 2-subgroup
$S$. Then, for all $n\in\Z$, the following induced map from the
group of coinvariants is an isomorphism:
\[
\incl_*: \UNil_n\infdec(\Z[S];\Z[S],\Z[S])_{F/S}
\xra{\quad} \UNil_n\infdec(\Z[F];\Z[F],\Z[F]).
\]
\end{thm}

\begin{cor}\label{Cor_OddOrder}
Suppose $F$ is a finite group of odd order.  Then the map induced by
the inclusion of the trivial subgroup is an isomorphism:
\[
\incl_*: \UNil_n\infdec(\Z; \Z, \Z) \xra{\quad}
\UNil_n\infdec(\Z[F]; \Z[F], \Z[F]).
\]
\end{cor}

These are summands of the surgery groups $L_n\infdec(\Z[\Gamma \x
D_\infty])$ for $\Gamma = 1, F$. The target $\UNil$ has a complete
set of invariants obtained first by the retraction $F \to 1$ to the
source, and then from the two-stage obstruction theory of
Connolly--Davis \cite{CD} (compare \cite{CR, BR}).

Historically, these types of finite groups have been studied for
classical $L$-theory.

\begin{rem}
Suppose $F$ is a finite group with a normal abelian Sylow
2-subgroup.
 Then the classical $L$-groups $L_*^h(\Z[F])$ are computed by Bak \cite{BakAbelian} and Wall
\cite[Cor. 2.4.3]{WallGroupRing}. The computation is extended to the
colimit $L_*\infdec(\Z[F])$ by Madsen--Rothenberg \cite{MR}.
Therefore, by considering one such group $F$ at a time, one can
determine
\[
L_*\infdec(\Z[F\x D_\infty]) \quad\text{and}\quad L_*\infdec(\Z[F][x]).
\]
\end{rem}

\subsection{Background}

An involution ${}^-$ on a  ring $R$ is an additive endomorphism that
reverses products and whose square is the identity. For each $n \geq
0$, the surgery obstruction group $L_n(R)$ is defined as cobordism
groups of $n$-dimensional quadratic Poincar\'e complexes over
$(R,{}^-)$ \cite{RanickiExact}.  Tensor product with the symmetric
Poincar\'e complex over $\Z$ of the complex projective plane $\CP^2$
induces a periodicity isomorphism $L_n(R) \to L_{n+4}(R)$, and so
the definition is extensible to all $n \in \Z$.  These abelian
groups fit into the surgery exact sequence \cite{Wall}, which can be
used to compute the set $\cS^s(X)/\hAut^s(X)$ of homeomorphism
classes of manifolds in the simple homotopy type of a closed
manifold $X$ of dimension $n > 4$.  Here, the ring $R =
\Z[\pi_1(X)]$ has involution given by inversion of group elements.

Let $\fB_-,\fB_+$ be $(R,R)$-bimodules with involution ${}^\wedge$
satisfying $\wh{r\cdot b\cdot s} = \ol{s} \cdot \wh{b} \cdot
\ol{r}$. The splitting obstruction group
$\UNil_{2k}^h(R;\fB_-,\fB_+)$ is defined as the Witt group of
interlocking quadratic forms over $\fB_\pm$ whose adjoints compose
to a nilpotent endomorphism. Product with the circle $S^1$ defines
the groups $\UNil_{2k-1}^h(R;\fB_-,\fB_+)$ \cite{CappellUnitary}.
These abelian groups satisfy 4-periodicity and are 2-primary.  In
the case of group rings of an amalgamated free product or
HNN-extension $G$, there is a split monomorphism $\UNil_n^h \longra
L_n^h(\Z[G])$ and $\UNil_n^h$ has exponent 8 \cite{Farrell}.  If the
manifold $X$ has a $\pi_1$-injective, two-sided submanifold $X_0$,
we obtain such a decomposition of fundamental groups.  A homotopy
equivalence $h: M \to X$ is \textbf{splittable} along $X_0$ if $h$
is homotopic to a map $h'$ transverse to $X_0$ and $(h')\inv(X_0)
\to X_0$ is a homotopy equivalence.  If $X_0$ has dimension $>4$,
Sylvain Cappell's nilpotent normal cobordism construction
\cite{CappellSplit} provides a bijection sending a homotopy
structure $[h] \in \cS(X)$ with vanishing Whitehead torsion in
$\Nil_0$ to a normally cobordant split solution $[h'] \in
\cS^\Split(X;X_0)$ and a splitting obstruction $[\Split(h)] \in
\UNil_{n+1}^h$ \cite{CappellFree}. The element $\Split(h)$ vanishes
if and only if $h$ is splittable along $X_0$.

Now, the source group of Corollary \ref{Cor_OddOrder} is computed
\cite[Theorem D]{BR} \cite[Theorem 1.10]{CD}.
\begin{thm}[Banagl, Connolly, Davis, Ranicki]
Let $n \in \Z$.  Then
\[
\UNil_n^h(\Z;\Z,\Z) = \begin{cases}0 & \text{if } n \equiv 0\pmod{4}\\
0 & \text{if } n \equiv 1\pmod{4}\\  x \Z[x]/2 & \text{if } n \equiv 2 \pmod{4}\\
\Z[x]/4 \x \bigoplus_3 \Z[x]/2 & \text{if } n \equiv 3
\pmod{4}.\end{cases}
\]
\end{thm}

The next corollary is immediate from the above discussion, Corollary
\ref{Cor_OddOrder}, and the semiperiodicity isomorphism
$\UNil_n^h(R;R^-,R^-) = \UNil_{n+2}^h(R;R,R)$ \cite{CappellUnitary}.

\begin{cor}
Let $Z$ be a closed, connected, oriented $n$-manifold with finite
fundamental group of odd order. For any $n+m
> 5$, define the codimension one submanifold
\[
X_0 := Z \x S^{m-1} \quad\text{of}\quad X := Z \x (\RP^m \# \RP^m).
\]
Here, the sphere $S^{m-1}$ defines the connected sum of the real
projective spaces.

Suppose $m$ is odd and $n+m \equiv 0,3 \pmod{4}$, or $m$ is even and
$n+m \equiv 1,2 \pmod{4}$. Then any simple homotopy equivalence $h:
M \to X$ is splittable along $X_0$, where $M$ is a closed manifold.
Otherwise, there are a countably infinite number of counterexamples
within any given cobordism class of $M$.\qed
\end{cor}

\subsection{Decorations and correspondences}

First, we describe a correspondence in algebraic $K$-theory. Let $R$
be a ring, and let $i \in \Z$. Using the ring map $\aug_0: R[x] \to
R$, Hyman Bass \cite{Bass1} defined a decomposition
\[
\wt{K}_i(R[x]) = \wt{K}_i(R) \x NK_i(R).
\]
Specializing to the case of an amalgamated product $R[D_\infty] =
R[C_2] *_R R[C_2]$ of rings, Friedhelm Waldhausen \cite{Waldhausen}
determined a decomposition
\[
K_i(R[D_\infty]) = \frac{K_i(R[C_2]) \x K_i(R[C_2])}{K_i(R)} \x \wt{\Nil}_{i-1}(R;R,R).
\]

These two $K$-theory $\Nil$-groups agree by the Davis--Khan--Ranicki
correspondence.

\begin{thm}[{\cite{DKR}}]
Let $R$ be a ring.  Then, for all $i \in \Z$, there is a natural
isomorphism
\[
\wt{\Nil}_i(R;R,R) \longra NK_{i+1}(R).
\]
\end{thm}

Next, we discuss the appropriate $K$-theory decorations in
$L$-theory \cite[Section 17]{RanickiLower}.

\begin{defn}
Let $R$ be a ring with involution.  For a given $i \leq 1$, let $S
\subseteq \wt{K}_i(R)$ and $T \subseteq NK_i(R)$ be $*$-invariant
subgroups. For each $n \in \Z$, define the \textbf{intermediate
lower, lower, and ultimate nilpotent $L$-groups} by
\begin{gather*}
NL_n^T(R) := \Ker\prn{\aug_0: L_n^{S\x T}(R[x]) \longra L_n^S(R)}\\
NL_n^{\gens{i}}(R) := NL_n^{NK_i(R)}(R)\\
NL_n^{\gens{-\infty}}(R) := \Ker\prn{\aug_0:
L_n^{\gens{-\infty}}(R[x]) \longra L_n^{\gens{-\infty}}(R)}.
\end{gather*}
\end{defn}
For typographical convenience, we shall abbreviate $NL_n(R) :=
NL_n\infdec(R)$.

\begin{defn}
Define the \textbf{simple, free, and projective nilpotent
$L$-groups} by
\begin{gather*}
NL_n^s(R) := NL_n^{0 \subseteq NK_1(R)}(R)\\
NL_n^h(R) := NL_n^{\gens{1}}(R) = NL_n^{0 \subseteq NK_0(R)}(R)\\
NL_n^p(R) := NL_n^{\gens{0}}(R).
\end{gather*}
\end{defn}

There are natural decompositions
\begin{gather*}
L_n^{S \x T}(R[x]) = L_n^S(R) \oplus NL_n^T(R)\\
L_n^{\gens{i}}(R[x]) = L_n^{\gens{i}}(R) \oplus
NL_n^{\gens{i}}(R)\\
L_n^{\gens{-\infty}}(R[x]) = L_n^{\gens{-\infty}}(R) \oplus
NL_n^{\gens{-\infty}}(R).
\end{gather*}

Using pairs of finitely generated free $R$-modules with additional
unitary structure, Sylvain Cappell \cite{CappellUnitary} defined the
$\UNil$-groups
\[
\UNil_n^h(R;R,R) = \UNil_n^{\wt{\Nil}_0(R;R,R)}(R;R,R)
\]
and showed that they fit into a decomposition
\[
L^h_n(R[D_\infty]) = \frac{L^h_n(R[C_2]) \x L^h_n(R[C_2])}{L^h_n(R)} \x \UNil_n^h(R;R,R).
\]

Following Cappell, we use a Shaneson-type sequence to define lower
decorations.

\begin{defn}
Write $\UNil_n^{\gens{1}}(R;R,R) := \UNil_n^h(R,R,R)$.  For every
$i\leq 0$, define
\[
\UNil_n^{\gens{i}}(R;R,R) :=
\frac{\UNil_n^{\gens{i+1}}(R[C_\infty];R[C_\infty],R[C_\infty])}{\UNil_n^{\gens{i+1}}(R;R,R)}.
\]
There is a forget-decoration map which fits into a Rothenberg-type
sequence
\[
\cdots \longra \UNil_n^{\gens{i+1}} \longra \UNil_n^{\gens{i}} \longra \wh{H}^n(C_2;\wt{\Nil}_{i-1})
 \xra{\en\bdry\en} \UNil_{n-1}^{\gens{i+1}} \longra \cdots.
\]
We define the \textbf{projective and ultimate unitary nilpotent
$L$-groups} by
\begin{gather*}
\UNil_n^p(R;R,R) := \UNil_n^{\gens{0}}(R;R,R)\\
\UNil_n^{\gens{-\infty}}(R;R,R) := \colim_{i \leq 1}
\UNil_n^{\gens{i}}(R;R,R).
\end{gather*}
\end{defn}

The earlier analogue to \cite{DKR} in $L$-theory was the
Connolly--Ranicki correspondence.

\begin{thm}[{\cite{CR}}]
Let $R$ be a ring with involution.
Then, for all $n \in \Z$, there is a natural isomorphism
\[
r^h: \UNil_n^h(R;R,R) \longra NL_n^h(R)
\]
which descends to natural isomorphisms
\begin{gather*}
r^p: \UNil_n^p(R;R,R) \longra NL_n^p(R)\\
r\infdec: \UNil_n\infdec(R;R,R) \longra NL_n\infdec(R).
\end{gather*}
\end{thm}

\subsection{Tools}

For the benefit of the reader, we list the technical tools for this
paper.
\begin{enumerate}
\item
The abelian groups $\UNil_*(R;R,R)$ vanish if 2 is a unit in $R$
(Theorem \ref{Thm_VanishNL}).
\item
Hyperelementary induction for the functor $NL_n(\Z[-])$ on
finite groups is concentrated at the prime $p=2$ (Theorem
\ref{Thm_DressInduction}).
\item
A cartesian square (a pullback-pushout diagram) of rings with
involution induces a Mayer--Vietoris exact sequence of
$NL$-groups (Theorem \ref{Thm_Pullback}).
\item
Chain bundles are used implicitly (Proposition
\ref{Prop_HomologicalReduction}) to prove that
$NL_*(\Z[\zeta_d])=0$ for all $d>1$ odd. Again, the emphasis
turns out to be at the prime 2, and we obtain that
$NL_*(\Z[\zeta_{e}])$ is a basic building block for all $e \geq
1$ a power of 2.
\end{enumerate}

\subsection{Motivation}

The Farrell--Jones isomorphism conjecture in $L$-theory \cite{FJIso}
states for any discrete group $\Gamma$ that $L_*(\Z[\Gamma])$ is
determined by $L_*(\Z[V])$ of all virtually cyclic subgroups $V$ of
$\Gamma$ together with certain homological information. That is, the
blocked assembly map is conjectured to be an isomorphism:
\[
H_n(B_{vc}(\Gamma); \bL_\bullet\infdec(\Z[-])) \longra L_n\infdec(\Z[\Gamma]).
\]
For computation of the source, there are spectral sequences
(Atiyah--Hirzebruch \cite{QuinnII} and $p$-chain Davis--L\"uck
\cite{DL}) which converge to the cosheaf $\bL$-homology of the
classifying space $B_{vc}(\Gamma)$ for $\Gamma$-actions with
virtually cyclic isotropy. A group $V$ is \textbf{virtually cyclic}
if it contains a cyclic subgroup of finite index. Equivalently:
\begin{enumerate}
\item[(I)]
$V$ is a finite group, or

\item[(II)]
$V$ is a group extension $1 \to F \to V \to C_\infty \to 1$ for
some finite group $F$, or

\item[(III)]
$V$ is a group extension $1 \to F \to V \to D_\infty \to 1$ for
some finite group $F$.
\end{enumerate}

The $L$-theory of type I, with various decorations, is determined
classically by Wall and others \cite{WallGroupRing} \cite{MR}. The
$L$-theory of type II is determined by a combination of the
$L$-theory of type I and the monodromy map $(1-\alpha_*)$ in the
Cappell--Ranicki--Shaneson--Wall sequence \cite[\S12B]{Wall}
\cite{RanickiLIII}. The groups $V$ of type III admit a decomposition
\[
V = V_- *_F V_+
\]
as an injective amalgam with
\[
[V_\pm : F] = 2.
\]
Sylvain Cappell developed a Mayer--Vietoris type sequence
\cite{CappellFree} that determines the groups $L_*^h(\Z[V])$ as an
extension of the $L$-groups of the type I groups $F, V_-, V_+$ and
of his splitting obstruction groups
\[
\UNil_*^h(\Z[F];\Z[V_-\setminus F], \Z[V_+\setminus F]).
\]
The recent computations \cite{CK, CR, CD, BR} of these
$\UNil$-groups for $F=1$ provide a starting point for our
determination of the $L$-groups $L_*(\Z[V])$ for certain classes of
type III virtually cyclic groups $V$.

\subsection{Outline of proof}
The main theorem (\ref{Thm_FiniteReduction}) and its corollary
(\ref{Cor_OddOrder}) are proven at the end of Section
\ref{Sec_FiniteReduction}.  The method is to apply 2-hyperelementary
induction (Section \ref{Sec_Hyperelementary}) and then to use
Mayer--Vietoris sequences for cartesian squares (Section
\ref{Sec_FiniteReduction}).  Odd order information is removed by
vanishing theorems (Section \ref{Sec_LocCompExc}) and by homological
analysis of cyclotomic number rings (Section
\ref{Sec_BasicReductions}).


\begin{ack}
The author would like to thank Jim Davis for frequent conversations
on the tools used in the classical $L$-theory of finite groups. This
material was adapted from the author's thesis
\cite{Khan_Dissertation} under his supervision at Indiana
University.  The referee and Andrew Ranicki both encouraged a close
attention to decorations.
\end{ack}

\section{Localization, completion, and excision for
$NL$}\label{Sec_LocCompExc}

In this section, we document a known vanishing result for $\UNil=NL$
(Theorem \ref{Thm_VanishNL}) and apply it to show that localization
and completion properties are concentrated at the prime 2 (Theorem
\ref{Thm_LocalComplete}). Also, we state the Mayer--Vietoris
sequence for a well-known cartesian square that allows us to pass
from the group ring of a group $G$ to the quotient $G/K$ by a finite
subgroup $K$ (Theorem \ref{Thm_Pullback}).

\begin{defn}
Let $j \in \Z$ and write $\epsilon := (-1)^j$. Define the
\textbf{Tate cohomology group}
\[
\wh{H}^j(A) := \wh{H}^j(C_2; A) = \frac{\set{a \in A \ST a =
\epsilon \ol{a}}}{\set{ b + \epsilon \ol{b} \ST b \in A }},
\]
where the group $C_2$ acts on the ring $A$ by its involution.
\end{defn}

We start with an elementary vanishing theorem (cf. Karoubi
\cite[Theorem 7]{Ojanguren}).  For other kinds of rings with
involution, it turns out that the vanishing of Tate cohomology
implies the vanishing of the nilpotent $L$-group (Proposition
\ref{Prop_HomologicalReduction}).

\begin{thm}\label{Thm_VanishNL}
Suppose 2 is a unit in a ring $A$ with involution. Let $j,n \in \Z$.
Then
\[
\wh{H}^j(A)=0 \quad\text{and}\quad N L_n(A) = 0.
\]
\end{thm}

The proof requires two algebraic lemmas.

\begin{lem}\label{Lem_BinomHalf}
For any $r \in \N$, the rational number $\binom{-1/2}{r}$ lies in
$\Z[\frac{1}{2}]$.
\end{lem}

\begin{pf}
This follows immediately from the observation
\[
\binom{-1/2}{r} = \binom{2r}{r} \prn{-\frac{1}{4}}^r,
\]
proven inductively, and from the fact that $\binom{2r}{r} \in \Z$.
\qed\end{pf}

\begin{lem}\label{Lem_SquareRoot}
Suppose $\rho$ is a nilpotent element of a ring $S$ that contains 2
as a unit. Then there exists a unit $V$ of $S$ commuting with $\rho$
such that $V^2 = (1+\rho)\inv$.
\end{lem}

\begin{pf}
Observe the nilpotence of $\rho \in S$ implies that $1+\rho \in
S^\times$, with inverse
\[
(1+\rho)\inv = \sum_{r=0}^\infty (-1)^r \rho^r.
\]
By Lemma \ref{Lem_BinomHalf} and the nilpotence of $\rho$, we can
define an element $V \in S$ that commutes with $\rho$:
\[
V := \sum_{r=0}^\infty \binom{-1/2}{r} \rho^r.
\]
Then the binomial theorem implies
\[
V^2 = (1+\rho)\inv.
\]
In particular, we obtain $V \in S^\times$. \qed\end{pf}

\begin{pf*}{PROOF of Theorem \ref{Thm_VanishNL}(1).}
Write $\epsilon := (-1)^j$. Suppose $a \in A$ is
$\epsilon$-symmetric. Define
\[
b := {\textstyle\frac{1}{2}} a \in A.
\]
Then note $a = \epsilon \ol{a}$ implies
\[
a = {\textstyle\frac{1}{2}}(a + a) = {\textstyle\frac{1}{2}} (a +
\epsilon \ol{a}) = b + \epsilon \ol{b}.
\]
Therefore every $\epsilon$-symmetric element of $A$ is also
$\epsilon$-even. Hence $\wh{H}^j(A)=0$. \qed\end{pf*}

\begin{pf*}{PROOF of Theorem \ref{Thm_VanishNL}(2).}
First consider the case $n=2k$ and write $\epsilon := (-1)^k$.
Recall that Higman linearization \cite[Lemma 4.2]{RanickiLIV}
involves stabilization by hyperbolic planes and zero-torsion
isomorphisms. Then we can represent
\[
\vartheta = [P[x], f_0 + x f_1] \in NL_{2k}^h(A)
\]
by a ``linear'' nonsingular $\epsilon$-quadratic form over $A[x]$
with null-augmentation to $A$. Here, $P$ is a finitely generated
free left $R$-module, and $f_0,f_1: P \longra P^* = \Hom_R(P,R)^t$
are left $R$-module morphisms. Recall, since 2 is a unit in $A[x]$,
that the $\epsilon$-symmetrization map is an isomorphism
\cite[Proposition 1.4.3]{RanickiExact}:
\[
L_{2k}^h(A[x]) \xra{\en 1+T_\epsilon \en} L^{2k}_h(A[x]).
\]
In fact, the quadratic refinement of a symmetric form is recovered,
uniquely up to skew $(-\epsilon)$-even morphisms, as one-half of the
$\epsilon$-symmetric morphism. So it is equivalent to show the
vanishing of the $\epsilon$-symmetric Witt class
\[
(1+T_\epsilon)(\vartheta) = [P[x], \lambda_0 + x \lambda_1].
\]
Here, for each $i=0,1$, the $\epsilon$-symmetrizations are defined
as
\[
\lambda_i := f_i + \epsilon f_i^*: P[x] \xra{\quad} P[x].
\]

After stabilization if necessary, there exists a lagrangian $P_0$ of
the $\epsilon$-symmetric form $(P,\lambda_0)$ over $A$, since
\[
\eval_0(\vartheta) = 0 \in L_{2k}^s(A).
\]
Since $\lambda_0 + x \lambda_1$ is invertible, we obtain a nilpotent
element of the ring $\End_{A[x]}(P[x])$:
\[
x \nu := \lambda_0\inv \circ x\lambda_1.
\]
Then by Lemma \ref{Lem_SquareRoot}, there exists
\[
V \in \End_{A[x]}(P[x])^\times
\]
commuting with $x\nu$ such that $V^2 = (1+ x\nu)\inv$. Observe that
$x\nu \in \End_{A[x]}(P[x])$ is a self-adjoint operator with respect
to the nonsingular form $(P[x],\lambda_0)$:
\[
(x \nu)^* \circ \lambda_0 = (\lambda_0\inv \circ x\lambda_1)^* \circ
\lambda_0 = x\lambda_1 \circ \lambda_0\inv \circ \lambda_0 = x
\lambda_1 = \lambda_0 \circ (x \nu).
\]
Hence the automorphism $V$ defined in Proof \ref{Lem_SquareRoot} is
also self-adjoint with respect to the symmetric form $\lambda_0$:
\[
V^* \circ \lambda_0 = \lambda_0 \circ V.
\]
Then note the pullback is $V^\%(P[x], \lambda_0 + x\lambda_1) =
(P[x],\lambda_0)$, since
\[
V^* \circ (\lambda_0 + x \lambda_1) \circ V = V^* \circ \lambda_0
\circ (1 + x \nu) \circ V = \lambda_0 \circ V \circ V \circ (1 + x
\nu) = \lambda_0.
\]
Hence the form $(P[x], \lambda_0 + x\lambda_1)$ is homotopy
equivalent to the symplectic form $(P[x],\lambda_0)$. That is,
$V(P_0[x])$ is a lagrangian for the $\epsilon$-symmetric form
$(P[x], \lambda_0 + x\lambda_1)$.  Therefore
\[
\vartheta = (1+T_\epsilon)\inv[P[x], \lambda_0 + x\lambda_1] = 0 \in
N L_{2k}^h(A).
\]

Thus we have shown $NL_{2k}^h(A)=0$.  By the fundamental theorem of
algebraic $L$-theory \cite{RanickiLIV}, we inductively obtain that
$NL_{2k}^{\gens{i}}(A)=0$ for all $i \leq 1$. Hence
$NL_{2k}\infdec(A)=0$.

Finally, the case $n=2k-1$ follows from the Ranicki--Shaneson
sequence \cite[Theorem 1.1]{RanickiLII} for ultimate decorations:
\[
0 \xra{\quad} NL_{2k}\infdec(A) \xra{\quad} NL_{2k}\infdec(A[C_\infty])
\xra{\quad} NL_{2k-1}\infdec(A) \xra{\quad} 0.
\]
\qed\end{pf*}

Next, we show that there are excision sequences for certain
cartesian squares.

\begin{lem}\label{Lem_MVsequences}
Suppose $\Phi$ is a cartesian square of rings with involution:
\[
\begin{diagram}
\node{A} \arrow{e} \arrow{s} \node{B} \arrow{s}\\
\node{C} \arrow{e} \node{D}
\end{diagram}
\]
such that $\Phi$ is a localization-completion square or that one of
the maps to $D$ is surjective. Then, for all $n \in \Z$, there are
Mayer--Vietoris exact sequences:
\begin{enumerate}
\item
of 2-periodic Tate cohomology groups
\[
\cdots \xra{\quad} \wh{H}^{n+1}(D) \xra{\en\bdry\en} \wh{H}^n(A)
 \xra{\quad} \wh{H}^n(B) \oplus \wh{H}^n(C) \xra{\quad} \wh{H}^n(D) \xra{\quad} \cdots
\]

\item
and of 4-periodic ultimate $NL$-groups
\[
\cdots \xra{\quad} NL_{n+1}(D) \xra{\en\bdry\en} NL_n(A)
\xra{\quad} NL_n(B) \oplus NL_n(C) \xra{\quad}
NL_n(D) \xra{\quad} \cdots.
\]
\end{enumerate}
\end{lem}

\begin{pf}
Given a finite group $\Gamma$, consider a short exact sequence of
$\Z[\Gamma]$-modules:
\[
0 \xra{\quad} M_0 \xra{\quad} M_- \oplus M_+ \xra{\quad} M
\xra{\quad} 0.
\]
Let $C$ be a contractible complex of f.g. free $\Gamma$-modules such
that
\[
\Cok(\bdry: C_1 \to C_0) = \Z,
\]
which is the trivial $\Z[\Gamma]$-module.  Recall, for any
coefficient $\Z[\Gamma]$-module $N$, the definition of Tate
cohomology:
\[
\wh{H}^j(\Gamma;N) = H^j(\Hom_{\Z[\Gamma]}(C,N)).
\]
Then we obtain the Bockstein sequence:
\[
\cdots \xra{\en\bdry\en} \wh{H}^j(\Gamma; M_0) \xra{\quad}
\wh{H}^j(\Gamma; M_-) \oplus \wh{H}^j(\Gamma; M_+) \xra{\quad}
\wh{H}^j(\Gamma; M) \xra{\en\bdry\en} \cdots.
\]
Therefore Part (1) is proven by substition: $\Gamma=C_2, M_0=A,
M_-=B, M_+=C, M=D$.

In order to prove Part (2), let $i \leq 0$. Consider the cartesian
square $\kappa$ of $*$-invariant subgroups
\cite[p.~498]{RanickiExact}:
\[
\kappa_i := \prn{\begin{diagram}
\node{I_i} \arrow{e} \arrow{s} \node{0} \arrow{s}\\
\node{0} \arrow{e} \node{0}
\end{diagram}}
\subseteq NK_i(\Phi)
\]
where
\[
I_i := \Ker\prn{NK_i(A) \longra NK_i(B) \oplus
NK_i(C)}.
\]
Recall, for any ring $R$, that the fundamental theorem of algebraic
$K$-theory \cite[XII]{Bass1} states
\[
K_i(R[x,x\inv]) = K_i(R) \x NK_i(R) \x NK_i(R) \x K_{i-1}(R).
\]
Consider the cartesian squares $\Phi[x], \Phi[x\inv], \Phi[x,x\inv]$
of the polynomial extensions and the decoration squares $\lambda_i,
\lambda_i, \mu_i$ defined by
\begin{gather*}
\lambda_i := \prn{\begin{diagram}
\node{J_i} \arrow{e} \arrow{s} \node{0} \arrow{s}\\
\node{0} \arrow{e} \node{0}
\end{diagram}}
\subseteq NNK_i(\Phi) \subseteq NK_i(\Phi[x]) \iso NK_i(\Phi[x\inv])\\
\mu_i := \kappa_i \x \lambda_i \x \lambda_i \x \kappa_{i-1}
\subseteq NK_i(\Phi[x,x\inv])
\end{gather*}
where
\[
J_i := \Ker\prn{NNK_i(A) \longra NNK_i(B) \oplus
NNK_i(C)}.
\]
Then, by the fundamental theorem of algebraic $L$-theory
($\ol{x}=x$) \cite[Corollary 4.4]{RanickiLIV}, we obtain a
decomposition of decorated triad $NL$-groups:
\[
NL_*^{\mu_i}(\Phi[x,x\inv]) = NL_*^{\kappa_i}(\Phi) \x NNL_*^{\lambda_i}(\Phi) \x NNL_*^{\lambda_i}(\Phi) \x NL_*^{\kappa_{i-1}}(\Phi).
\]
More generally, write $X_i := \set{x_1, \ldots, x_{1-i}}$. Thus, we
inductively obtain that $NL_*^{\kappa_{i-1}}(\Phi)$ is a summand of
$NL_*^{\omega_0}(\Phi[X_i, X_i\inv])$ for a certain cartesian square
$\omega_0 \subseteq NK_0(\Phi[X_i, X_i\inv])$ of $*$-invariant
subgroups. Observe that $\Phi[X_i, X_i\inv]$ is a cartesian square
of rings with involution satisfying the same stated condition as
$\Phi$. Then, by Ranicki's excision theorem for intermediate
projective decorations \cite[Proposition 6.3.1]{RanickiExact}, we
obtain $NL_*^{\kappa_0}(\Phi)=0$ and
$NL_*^{\omega_0}(\Phi[X_i,X_i\inv]) = 0$. Hence
$NL_*^{\kappa_i}(\Phi)=0$ for all $i \leq 0$.  In other words, we
obtain a Mayer--Vietoris exact sequence
\[
\cdots \xra{\quad} NL_{n+1}^{\gens{i+1}}(D) \xra{\en\bdry\en} NL_n^{I_i}(A)
\xra{\quad} NL_n^{\gens{i+1}}(B) \oplus NL_n^{\gens{i+1}}(C) \xra{\quad} NL_n^{\gens{i+1}}(D) \xra{\quad}
\cdots.
\]
Now, observe the intermediate and full decorations intertwine by
forget-decoration maps:
\[
NL_n^{I_0}(A) \longra NL_n^{NK_0(A)}(A) \longra NL_n^{I_{-1}}(A) \longra NL_n^{NK_{-1}(A)}(A) \longra \cdots.
\]
This cofinality implies that the induced map between the direct
limits is an isomorphism:
\[
\colim_{i \leq 0} NL_n^{I_i}(A) = \colim_{i\leq 0} NL_n^{NK_i(A)}(A)
= NL_n\infdec(A).
\]
Therefore we obtain the desired exact sequence of
$NL_*\infdec$-groups as a direct limit. \qed
\end{pf}

The following basic isomorphism shall be useful for ring
decompositions of group rings.

\begin{thm}\label{Thm_LocalComplete}
Let $A$ be a ring with involution, and let $n \in \Z$.

\begin{enumerate}
\item
For all $N$ odd, the following induced maps are isomorphisms:
\[
\wh{H}^n(A) \xra{\quad} \wh{H}^n(A[\textstyle\frac{1}{N}])
\quad\text{and}\quad NL_n(A) \xra{\quad}
NL_n(A[\textstyle\frac{1}{N}]).
 \]

\item
The following induced maps are isomorphisms:
\[
\wh{H}^n(A) \xra{\quad} \wh{H}^n(\wh{A}_{(2)}) \quad\text{and}\quad
NL_n(A) \xra{\quad} NL_n(\wh{A}_{(2)}).
\]
\end{enumerate}
\end{thm}

\begin{pf}
Let $r > 0$. Consider the localization-completion cartesian square
$\Phi$ of rings with involution \cite[pp. 197--8]{RanickiExact}:
\[
\Phi := \prn{
\begin{diagram}
\node{A} \arrow{e} \arrow{s} \node{\wh{A}_{(r)}} \arrow{s}\\
\node{A[\textstyle\frac{1}{r}]} \arrow{e}
\node{\wh{A}_{(r)}[\textstyle\frac{1}{r}]}
\end{diagram}
}.
\]
Note that $2$ is a unit in the localization
$A[\textstyle\frac{1}{2}]$ and in the completion $\wh{A}_{(N)}$.
Then, by Theorem \ref{Thm_VanishNL}, we obtain:
\begin{gather*}
\wh{H}^*(A[\textstyle\frac{1}{2}]) = \wh{H}^*(\wh{A}_{(2)}[\textstyle\frac{1}{2}])
 = 0 = \wh{H}^*(\wh{A}_{(N)}) = \wh{H}^*(\wh{A}_{(N)}[\textstyle\frac{1}{2}])\\
NL_*(A[\textstyle\frac{1}{2}]) =
NL_*(\wh{A}_{(2)}[\textstyle\frac{1}{2}])
 = 0 = NL_*(\wh{A}_{(N)}) = NL_*(\wh{A}_{(N)}[\textstyle\frac{1}{2}]).
\end{gather*}
Therefore we are done by the Mayer--Vietoris exact sequences of
Lemma \ref{Lem_MVsequences}. \qed
\end{pf}

Now we establish useful Mayer--Vietoris sequences for certain normal
subgroups.

\begin{thm}\label{Thm_Pullback}
Let $R$ be a ring with involution of characteristic zero. Let $K$ be
a finite normal subgroup of any group $G$. Its \textbf{norm} is
defined as
\[
\Sigma_K := \sum_{g \in K} g \in R[G].
\]
\begin{enumerate}
\item
There is a cartesian square of rings with involution:
\[\begin{diagram}
\node{R[G]} \arrow{e} \arrow{s} \node{R[G]/\Sigma_K} \arrow{s}\\
\node{R[G/K]} \arrow{e} \node{\textstyle\frac{R}{\abs{K}}[G/K].}
\end{diagram}\]

\item
There is a Mayer--Vietoris exact sequence of 2-periodic Tate
cohomology groups:
\begin{multline*}
\cdots \xra{\quad} \wh{H}^n(R[G]) \xra{\quad} \wh{H}^n(R[G/K])
\oplus \wh{H}^n(R[G]/\Sigma_K) \\ \xra{\quad}
\wh{H}^n(\textstyle\frac{R}{\abs{K}}[G/K]) \xra{\en\bdry\en}
\wh{H}^{n-1}(R[G]) \xra{\quad} \cdots.
\end{multline*}

\item
There is a Mayer--Vietoris exact sequence of 4-periodic ultimate
$NL$-groups:
\begin{multline*}
\cdots \xra{\quad} NL_n(R[G]) \xra{\quad} NL_n(R[G/K]) \oplus
NL_n(R[G]/\Sigma_K)
\\ \xra{\quad} NL_n(\textstyle\frac{R}{\abs{K}}[G/K])
\xra{\en\bdry\en} NL_{n-1}(R[G]) \xra{\quad} \cdots.
\end{multline*}
\end{enumerate}
\end{thm}

\begin{pf}
It is straightforward to check that the diagram of Part (1) is
commutative and a pushout-pullback square (that is, a cartesian
square). Therefore we are done by the Mayer--Vietoris exact
sequences of Lemma \ref{Lem_MVsequences}. \qed\end{pf}

\section{Hyperelementary induction for $NL$}\label{Sec_Hyperelementary}

In this section, we use categorical language to document the Mackey
Subgroup Property for $\UNil=NL$ (Lemma \ref{Lem_Mackey}) and show
that only 2-hyperelementary subgroups of a finite group $F$ are
required for induction (Theorem \ref{Thm_DressInduction}).  Also, we
decompose $NL$ for a 2-hyperelementary group $H$ into components
involving cyclotomic number rings $\Z[\zeta_d]$ (Theorem
\ref{Thm_HyperDecomp}).

For any prime $p$, recall that a finite group $H$ is
\textbf{$p$-hyperelementary} if it is a group extension
\[
1 \xra{\quad} C_N \xra{\quad} H \xra{\quad} P \xra{\quad} 1
\]
with $P$ a finite $p$-group and $N$ coprime to $p$; the extension is
necessarily split.

For every group $G$ and prime $p$, define a category $\cH_p(G)$ as
follows. Its objects $H$ are all the $p$-hyperelementary subgroups
of $G$, and its morphisms are defined as conjugate-inclusions:
\[
\vphi_{H,g,H'} : H \xra{\quad} H';\qquad x \longmapsto g x g\inv,
\]
for all possible $p$-hyperelementary subgroups $H,H'$ and elements
$g$ of $G$.  It is a subcategory of the category $\FINITEGROUPS$ of
finite groups and monomorphisms. For any normal subgroup $S$ of $G$,
the category $\cH_p(G)$ contains a full subcategory $\cH_p(G)\cap S$
whose objects are all the $p$-hyperelementary subgroups of $S$. The
inclusion map also induces a functor
\[
\incl_*: \cH_p(S) \xra{\quad} \cH_p(G)\cap S,
\]
which is bijective on object sets and injective on morphism sets.

\begin{defn}\label{Def_MackeyNL}
Fix $n\in\Z$. Define a pair $\fN = (\fN^*, \fN_*)$ of functors by
\begin{gather*}
\fN_* : \FINITEGROUPS \xra{\quad} \ABELIANGROUPS\\
\fN_*(F) := NL_n\infdec(\Z[F])\\
\fN_*\prn{F\xra{\vphi}G} := \prn{ \fN(F) \xra{\vphi_*} \fN(\vphi F)
\xra{\incl_*} \fN(G) }
\end{gather*}
and
\begin{gather*}
\fN^* : \FINITEGROUPS^\op \xra{\quad} \ABELIANGROUPS\\
\fN^*(F) := NL_n\infdec(\Z[F])\\
\fN^*\prn{F\xra{\vphi}G} := \prn{ \fN(G) \xra{\incl^*} \fN(\vphi F)
\xra{\vphi_*\inv} \fN(F) }.
\end{gather*}
The morphism $\incl^*$ is the transfer map from $G$ to the finite
index subgroup $\vphi F$ \cite[Remark 21.7]{RanickiTOP}.  We write
$\fN(F)$ for the common value $\fN^*(F)=\fN_*(F)$.
\end{defn}

We shall reduce to colimits and coinvariants are consequences of
Dress induction \cite{Dress} and Farrell's exponent theorem
\cite{Farrell}.  Our theorem is proven over the next few pages.

\begin{thm}\label{Thm_DressInduction}
Let $F$ be a finite group and $S$ a normal subgroup.
\begin{enumerate}
\item
The following induced map from the direct limit is an isomorphism:
\[
\incl_*: \colim_{\cH_2(F)} \fN \xra{\quad} \fN(F).
\]

\item
The following map, from the group of coinvariants, is an induced
isomorphism:
\[
\incl_*: \prn{\colim_{\cH_2(S)} \fN}_{F/S} \xra{\quad}
\colim_{\cH_2(F)\cap S} \fN.
\]
\end{enumerate}
\end{thm}

\begin{defn}
Let $F$ be a finite group and $R$ a Dedekind domain.  A
\textbf{lattice} $(M,\lambda)$ consists of an $R$-projective,
finitely generated left $R[F]$-module $M$ and symmetric $R$-bilinear
map $\lambda: M \x M \to R$ that induces a left $R[F]$-module
isomorphism $\Ad(\lambda): M \to \Hom_R(M,R)^t$.  The
\textbf{equivariant Witt ring} $GW_0(F,R)$ is the Witt ring of
$F$-lattices over $R$, with the operations of direct sum and tensor
product over $R$.
\end{defn}

\begin{lem}\label{Lem_Mackey}
The above $\fN$ is a Mackey functor and a module over the Green ring
functor
\[
GW_0(-,\Z) : \FINITEGROUPS \xra{\quad} \COMMUTATIVERINGS.
\]
\end{lem}

\begin{pf}
Observe that $\fN$ transforms inner automorphisms into the identity,
by Taylor's Lemma \cite[Cor. 1.1]{Taylor}, since we are assuming
that all our groups are equipped with trivial orientation character.
Also, for any isomorphism $\vphi: F \to G = \vphi(F)$, by Definition
\ref{Def_MackeyNL}, we have
\[
\fN^*(\vphi) = \prn{\fN(G) \xra{\vphi_*\inv} \fN(F)} = \prn{\fN(F)
\xra{\vphi_*} \fN(G)}\inv = \fN_*(\vphi)\inv.
\]

Next we document the Mackey subgroup property (compare \cite[Theorem
4.1]{BakEvenOdd}). However, we shall more generally do so for the
analogously defined (\ref{Def_MackeyNL}) quadratic $L$-theory pair
of functors
\[
\fL = (\fL^*,\fL_*) = L_n\infdec(R[-]) : \FINITEGROUPS
\xra{\quad} \ABELIANGROUPS
\]
for any ring $R$ with involution. Let $H,K$ be subgroups of a finite
group $F$. Then we must show that the ``double coset formula''
holds, i.e. the following diagram commutes (we suppress labels for
the inclusions):
\[
\begin{diagram}
\node{\fL(H)} \arrow{e,t}{\fL_*} \arrow{s,t}{\fL^*} \node{\fL(F)}
\arrow{e,t}{\fL^*}
\node{\fL(K).}\\
\node{\bigoplus_{KaH \in K\backslash G/H} \fL(a\inv K a \cap H)}
\arrow[2]{ee,t}{(\Conj^a_*)} \node{} \node{\bigoplus_{KaH \in
K\backslash F/ H} \fL(K \cap a H a\inv)} \arrow{n,t}{\fL_*}
\end{diagram}
\]

Let $P$ be an arbitrary left $R[H]$-module.  Denote the inclusions
$i: H \hookra F$ and $j: K \hookra F$.  For all double cosets $K a
H$ with $a \in F$, denote inclusions
\[
i_a : K \cap a H a\inv \hookra a H a\inv \quad\text{and}\quad j_a: a
H a\inv \hookra K.
\]
The Mackey Subgroup Theorem \cite[Theorem 44.2]{CurtisReiner}
states, as an internal sum of $R[K]$-modules, that
\[
j^* i_* (P) = \bigoplus_{KaH \in K\backslash F/H} j_{a*} i_a^* (a\xo
P).
\]
Observe
\[
\Conj^a_*(P) = a \xo P \subseteq R[F] \xo_{R[H]} P,
\] where the
$R[aHa\inv]$-module structure on the $R$-submodule $a\xo P$ is given
by
\[
(a x a\inv) \cdot (a \xo p) = a \xo x p.
\]
Denote an inclusion
\[
k_a : a\inv K a \cap H \hookra H.
\]
Since $(\Conj^a)^* = (\Conj^a_*)\inv$ and the following diagram
commutes:
\[
\begin{diagram}
\node{a\inv K a \cap H} \arrow{e,t}{k_a} \arrow{s,tb}{\Conj^a}{\iso} \node{H} \arrow{s,tb}{\Conj^a}{\iso}\\
\node{K \cap a H a\inv} \arrow{e,t}{i_a} \node{a H a\inv}
\end{diagram}
\]
we obtain
\[
i_a^*(a\xo P) = i_a^* \Conj^a_*(P) = \Conj^a_* k_a^*(P).
\]
Hence the Mackey Subgroup Theorem is equivalent to the formula
\[
j^* i_*(P) = \bigoplus_{KaH \in K\backslash F/H} j_{a*} \Conj^a_*
k_a^*(P),
\]
and is functorial in left $R[H]$-modules $P$. Now consider its dual
module
\[
P^* := \Hom_{R[H]}(P, R[H])^t.
\]
There is a functorial $R[K]$-module morphism, which is an
isomorphism if $P$ is finitely generated projective:
\[
\Phi_a: j_a^* \Conj^a_* k_a^* (P^*) \xra{\quad} j_a^* \Conj^a_*
k_a^*(P)^*;\quad r \xo a \xo f \longmapsto (s \xo a \xo x \mapsto
s\, k_a^! f(x)\: \ol{r}).
\]
Here, the trace
\[
k_a^!: R[H] \to R[a\inv K a \cap H]
\]
is defined as projection onto the trivial coset (see \cite[Example
5.15]{HRT}). Thus for all f.g. projective $R[H]$-modules $P$, we
obtain a functorial isomorphism, which respects the above double
coset decomposition:
\[
\Phi := \prod_{KaH} \Phi_a: j^* i_*(P^*) \xra{\quad} j^* i_*(P)^*.
\]
Then there is a commutative diagram of algebraic bordism categories
and their functors \cite[\S3]{RanickiTOP}:
\[
\begin{diagram}
\node{\Lambda(R[H])} \arrow{e,t}{i_*} \arrow{s,t}{\prod k_a^*}
\node{\Lambda(R[G])} \arrow{e,t}{j^*}
\node{\Lambda(R[K])}\\
\node{\prod_{KaH} \Lambda(R[a\inv K a \cap H])} \arrow[2]{e,t}{\prod
\Conj^a_*} \node{} \node{\prod_{KaH} \Lambda(R[K \cap a H a\inv])}
\arrow{n,t}{\prod j_{a*}}
\end{diagram}
\]
with $\gens{-\infty}$ decorations.  So the desired commutative
diagram is induced \cite[Proposition 3.8]{RanickiTOP} on the level
of $L_*\infdec$-groups. Therefore $\fL$ (resp. $\fN$) is a Mackey
functor. The module structure on $\fL$ (resp. $\fN$) over the Green
ring functor
\[
GW_0(-,\Z)
\]
is defined (see \cite[p. 1452]{BakEvenOdd}, resp. \cite[p.
306]{Farrell}) using the diagonal $F$-action:
\[
GW_0(F,\Z) \x \fL(F) \xra{\quad} \fL(F); \qquad \prn{[M,\lambda],
[C,\psi]} \mapsto [M \xo_{\Z} C, \Ad(\lambda) \xo \psi].
\]
This verifies all assertions for $\fN$. \qed\end{pf}

\begin{pf*}{PROOF of Theorem \ref{Thm_DressInduction}(1).}
Since Lemma \ref{Lem_Mackey} shows that Dress Induction
\cite[Theorem 1]{Dress} is applicable in its covariant form
\cite[Theorem 11.1]{Oliver}, the functor $\Z_{(2)} \xo \fN$ is
$\cH_2$-computable. That is, the following induced map is an
isomorphism:
\[
\incl_*: \colim_{\cH_2(F)} \Z_{(2)} \xo \fN \xra{\quad} \Z_{(2)} \xo
\fN(F).
 \]
But \cite[Theorem 1.3]{Farrell} states for all groups $G$ that
$\fN(G)$ has exponent 8. Hence the prime 2 localization map
\[
\fN(G) \xra{\quad} \Z_{(2)} \xo \fN(G)
\]
is an isomorphism. The result follows immediately. \qed\end{pf*}

\begin{pf*}{PROOF of Theorem \ref{Thm_DressInduction}(2).}
For existence and surjectivity of the map, it suffices to show that
the following commutative diagram exists:
\[
\begin{diagram}
\node{\colim_{\cH_2(S)} \fN} \arrow{s,t,A}{\quotient} \arrow{e,t,A}{\incl_*} \node{\colim_{\cH_2(F)\cap S} \fN.}\\
\node{\prn{\colim_{\cH_2(S)} \fN}_{F/S}} \arrow{ne,t,..}{\incl_*}
\end{diagram}
\]
The group $F$ has a covariant action on the category $\cH_2(S)$
defined by pushforward along conjugation:
\[
\Conj: F \xra{\quad} \Aut(\cH_2(S)).
\]
Recall that $\fN$ transforms inner automorphisms into the identity,
by Taylor's Lemma \cite[Cor. 1.1]{Taylor}, since we are assuming
that all our groups are equipped with trivial orientation character.
Then the quotient group $F/S$ has an induced action on the colimit.
The \textbf{group of coinvariants} is defined by
\[
\prn{\colim_{\cH_2(S)} \fN}_{F/S} := \prn{\colim_{\cH_2(S)} \fN}
\bigg/ \bigg\langle{x_H - \Conj^g_*(x_H) \;\bigg|\; x_H \in \fN(H)
\text{ and } gS \in F/S}\bigg\rangle.
\]
Recall that the direct limit of a functor is defined in this case by
\[
\colim_{\cH_2(S)} \fN := \prn{\prod_{H \in \Ob\,\cH_2(S)} \fN(H)}
\bigg/ \bigg\langle{ x_H - \fN_*(\vphi)(x_H) \;\bigg|\; \vphi \in
\Mor\,\cH_2(S) }\bigg\rangle,
\]
and similarly over the finite category $\cH_2(F) \cap S$. This
explains the terms in the above diagram.

In order to show that the induced map exists, let
\[
z := x_H - \Conj^g_*(x_H) \in \prod \fN(H)
\]
represent a generator of the kernel of the quotient map, where $x_H
\in \fN(H)$ and $g \in F$. But note
\[
z = x_H - \fN_*(\vphi_{H,g,gHg\inv})(x_H).
\]
Hence it maps to zero in the direct limit over $\cH_2(F)\cap S$.
Thus the desired map exists and is surjective.

In order to show that the induced map is injective, suppose $[x_H]$
is an equivalence class in the coinvariants which maps to zero in
the direct limit over $\cH_2(F)\cap S$. Then there exists an
expression
\[
(x_H) = \sum_{i=1}^r (w_i - \fN_*(\vphi_i)(w_i)) \in \prod_{H \in
\Ob\,\cH_2(S)} \fN(H)
\]
for some
\[
H_1,\ldots, H_r \in \Ob\,\cH_2(S) \quad\text{and}\quad w_i \in
\fN(H_i) \quad\text{and}\quad \vphi_i \in \Mor\,\cH_2(F)\cap S.
\]
But each monomorphism $\vphi_i = \vphi_{H_i, g_i, H_i'}$ admits a
factorization
\[
\vphi_i = \vphi_i' \circ \Conj^{g_i}_*
\]
into an isomorphism $\Conj^{g_i}_*$ and an inclusion
\[
\vphi_i' := \vphi_{g_i H_i g_i\inv, 1, H'} \in \Mor\,\cH_2(S).
\]
Then note
\[
[w_i - \fN_*(\vphi_i)(w_i)] = [w_i - \Conj^{g_i}_*(w_i)] + [v_i -
\fN_*(\vphi_i')(v_i)] = 0 \in \prn{\colim_{\cH_2(S)} \fN}_{F/S}
\]
where
\[
v_i := \Conj^{g_i}_*(w_i).
\]
Hence $[x_H] = 0$ in the coinvariants. Thus the desired map
$\incl_*$ is injective. \qed\end{pf*}


Therefore, we are reduced to the computation of
\[
\fN(H) = NL_*\infdec(\Z[H])
\]
for all 2-hyperelementary groups $H$.  The next theorem reduces it
further.

\begin{thm}\label{Thm_HyperDecomp}
Suppose $H$ is a 2-hyperelementary group:
\[
H = C_N \rtimes_\tau P
\]
where $N$ is odd and $P$ is a finite 2-group. Consider the ring $R
:= \Z[\frac{1}{N}]$. Then for all $n \in \Z$, there is an induced
isomorphism
\[
NL_n(\Z[H]) \xra{\quad} NL_n\prn{(\bigoplus_{d|N} R[\zeta_d])
\circ_{\tau'} P},
\]
where the action $\tau'$ is induced by $\tau$.

Moreover if $\tau$ is trivial, then it can be lifted to an induced
isomorphism
\[
NL_n(\Z[C_N \x P]) \xra{\quad} \bigoplus_{d|N}
NL_n\prn{\Z[\zeta_d][P]}.
\]
\end{thm}

Suppose $A$ is a ring, $P$ is a group, and $\tau: P \to \Aut(A)$ is
a homomorphism.  So $P$ acts on $A$ by ring automorphisms. Then $A
\circ_\tau P$ denotes the \emph{twisted group ring} of $P$ with
coefficients from $A$.

\begin{pf}
For each divisor $d$ of $N$, let $\rho_d$ be the cyclotomic
$\Q$-representation of $C_N$ defined by $\rho_d(T) := \zeta_d$. The
$\rho_d$ represent all the distinct isomorphism classes of
irreducible $\Q$-representations of the group $C_N$.  An elementary
argument using the chinese remainder theorem and the trace of right
multiplication shows that there exists a restricted isomorphism of
$R$-algebras with involution:
\[
\rho = \bigoplus_{d|N} \rho_d: R[C_N] \xra{\quad} \bigoplus_{d|N}
R[\zeta_d].
\]

Therefore, using Theorem \ref{Thm_LocalComplete}, we obtain a
composition of isomorphisms:
\[
NL_n(\Z[C_N] \circ_\tau P) \xra{\quad} NL_n\prn{R[C_N] \circ_\tau
P} \xra{\en\rho_*\en} NL_n\prn{(\bigoplus_{d|N} R[\zeta_d])
\circ_{\tau'} P}.
\]
The assertion for $\tau$ trivial follows from Theorem
\ref{Thm_LocalComplete} and the additivity of $L_*\infdec$ (hence
$NL_*\infdec$) for finite products of rings with involution (cf.
\cite[Cor. 5.13]{HRT}). \qed\end{pf}


\section{Basic reductions}\label{Sec_BasicReductions}

Continuing on, we establish four reductions (orientable,
hyperelementary, nilpotent, homological) that we shall use to prove
Theorem \ref{Thm_FiniteReduction}.

\subsection{Orientable reduction}

\begin{prop}\label{Prop_OrientableReduction}
Suppose $R$ is a ring with involution and $G$ is a group with
trivial orientation character.  Then there is a natural
decomposition
\[
NL_n(R[G]) = NL_n(R) \oplus \wt{NL}_n(R[G]),
\]
where the \textbf{reduced $L$-group} is defined by
\[
\wt{NL}_n(R[G]) := \Ker\prn{\aug_1: NL_n(R[G]) \to
NL_n(R[1])}.
\]
\end{prop}

\begin{pf}
The covariant morphism on $NL_*(R[-])$-groups induced by the map
$\incl: 1 \to G$ is a monomorphism split by the morphism of
$NL_*(R[-])$-groups induced by the augmentation $\aug_1: G \to 1$ of
groups with orientation character. \qed\end{pf}

\subsection{Hyperelementary reduction}

For simplicity, from Section \ref{Sec_Hyperelementary}, we shall
continue the notation for fixed $n \in \Z$:
\[
\fN(G) := NL_n\infdec(\Z[G]).
\]

\begin{thm}\label{Thm_HyperelementaryReduction}
Let $F$ be a finite group and $S$ a normal subgroup. Suppose for all
2-hyperelementary subgroups $H$ of $F$ that the following induced
map is an isomorphism:
\[
\incl_*: \fN(H \cap S) \xra{\quad} \fN(H).
\]
Then the following induced map, from the group of coinvariants, is
an isomorphism:
\[
\incl_*: \fN(S)_{F/S} \xra{\quad} \fN(F).
\]
\end{thm}

\begin{pf}
Observe that the following diagram commutes:
\[
\begin{diagram}
\node{\colim_{\cH_2(F)\cap S} \fN} \arrow{se,t}{\incl_*}\\
\node{\prn{\colim_{\cH_2(S)}\fN}_{F/S}}
\arrow{n,t}{\incl_*} \arrow{e,t}{\incl_*} \arrow{s,t}{\incl_*} \node{\colim_{\cH_2(F)} \fN} \arrow{s,t}{\incl_*}\\
\node{\fN(S)_{F/S}} \arrow{e,t}{\incl_*} \node{\fN(F).}
\end{diagram}
\]
The vertical maps are isomorphisms by Hyperelementary Induction
(\ref{Thm_DressInduction}).  It follows from the hypothesis that the
diagonal map is an isomorphism.  Therefore the bottom map is an
isomorphism. \qed\end{pf}

\subsection{Nilpotent reduction}

The following is a specialization of Wall's theorem for complete
semilocal rings \cite[Theorem 6]{WallCompleteSemilocalRings}, which
was applied extensively in \cite{WallGroupRing}.  In the classical
$L$-theory of finite groups, Wall's theorem was applied to the
Jacobson radical \cite[\S 3]{WallCompleteSemilocalRings} of the
2-adic integral group ring of finite 2-groups \cite[\S
5.2]{WallGroupRing}.

In our case, a theorem of Amitsur \cite[Theorem 1]{Amitsur} states
that the Jacobson radical of $R[x]$ for any ring $R$ is a two-sided
ideal $N[x]$, where $N$ is a nil ideal of $R$ containing the locally
nilpotent radical.  Recall for left artinian rings $R$ that its
locally nilpotent radical, nilradical, and Jacobson radical all
coincide. In our applications, we limit ourselves to rings of
nonzero characteristic.

\begin{prop}\label{Prop_NilpotentReduction}
Let $R$ be a ring with involution.

\begin{enumerate}
\item
Suppose that $I$ is a nilpotent, involution-invariant, two-sided
ideal of $R$. Then for all $n\in\Z$, the map induced by the quotient
map $\pi: R \to R/I$ is an isomorphism:
\[
\pi_*: NL_n^h(R) \xra{\quad} NL_n^h(R/I).
\]

\item
Suppose for some prime $p$ that $\F$ is a finite field of
characteristic $p$ and $P$ is a finite $p$-group. Then for all $n\in
\Z$, the following induced map is an isomorphism:
\[
\aug_1: NL_n^h(\F[P] \xo_{\Z} R) \xra{\quad} NL_n^h(\F \xo_{\Z} R).
\]
\end{enumerate}
\end{prop}

\begin{pf}
For Part (1), observe that $I[x]$ is nilpotent implies that the map
\[
R[x] \longra \wh{R[x]}_{I[x]},
\]
to the $I[x]$-adic completion, is an isomorphism of rings with
involution. Then, since $\wt{K}_1(R[x]) = (\pi_*)\inv
\wt{K}_1(R[x]/I[x])$, by \cite[Theorem
6]{WallCompleteSemilocalRings}, we have induced isomorphisms:
\begin{gather*}
\pi_*: L_n^h(R[x]) \longra L_n^h(R[x]/I[x])\\
\pi_*: L_n^h(R) \longra L_n^h(R/I).
\end{gather*}
Therefore we obtain that $\pi_*: NL_n^h(R) \to NL_n^h(R/I)$ is an
isomorphism.

For Part (2), observe that the involution-invariant, two-sided ideal
\[
J := \prn{\set{g-1 \ST g \in P}}
\]
of $\F[P]$ is its Jacobson radical. Since $\F[P]$ is finite hence
left artinian, we must have that $J$ is nilpotent. So we are done by
Part (1) using $I=J\xo 1_R$. \qed\end{pf}

\subsection{Homological reduction}\label{Subsec_HomRed}

The following little observation shall be applied to algebraic
number ring $\cO$ with a Galois involution in the next section.

\begin{prop}\label{Prop_HomologicalReduction}
Let $R$ be a Dedekind domain with involution. Suppose the Tate
cohomology groups vanish: $\wh{H}^*(R)=0$. Then the nilpotent
$L$-groups vanish: $NL_*^h(R)=0$.
\end{prop}

\begin{pf}
Let $n \in \Z$.  By \cite[Proposition 20]{CR}, there is an
isomorphism
\[
NL_n^h(R) \iso NQ_n^h(R) := \Ker\prn{ Q_n(B^{R[x]},\beta^{R[x]}) \longra Q_n(B^R, \beta^R) }.
\]
The pair $(B^A,\beta^A)$ is the so-called universal chain bundle of
Michael S.~Weiss \cite[2.A.4]{WeissKervaireInvariant}.  These are
constructed for any ring $A$ with involution, by first selecting
free $A$-modules $F_k$ and $A$-module epimorphisms $F_k \to
\wh{H}^k(A)$ for every $k \in \Z$.  Here, the left $A$-module
structure is defined by
\[
A \x \wh{H}^k(A) \longra \wh{H}^k(A);\quad
(\lambda, [a]) \longmapsto [\lambda a \ol{\lambda}].
\]
However, we have assumed that $\wh{H}^*(R)=0$, which implies
$\wh{H}^*(R[x])=0$ by direct calculation.  Here, the involution on
$R[x]$ is extended from $R$ by $\ol{x}=x$. Then we can take $F_k=0$
in Weiss's construction for both rings $R$ and $R[x]$ with
involution. So the $\Z$-graded free module chain complexes vanish:
$B^{R}=0$ and $B^{R[x]} = 0$. Therefore $NL_n^h(R)=0$. \qed\end{pf}

\section{Finite groups with normal abelian Sylow 2-subgroup}\label{Sec_FiniteReduction}

The goal of this section is to prove Theorem
\ref{Thm_FiniteReduction}.  Observe that any finite group $F$
satisfying the hypothesis is a group extension
\[
1 \longra S \longra F \longra E \longra 1
\]
for a unique finite, abelian 2-group $S$ and odd order group $E$.
Since $H^2(E;S) = 0$ by a transfer argument \cite[Cor. 3.13]{Brown},
$F$ must be of the form $F = S \rtimes E$.

Many techniques \cite{WallGroupRing} \cite{HM} used to compute the
quadratic $L$-theory of finite groups, namely: hyperelementary
induction, the Mayer--Vietoris sequence for cartesian squares,
nilradical quotients, maximal involuted orders, and Morita
equivalence, along with our new technique of homological reduction
(Section \ref{Subsec_HomRed}), are employed in combination to
determine the quadratic $NL$-theory of certain finite groups, up to
extension issues.

The first lemma is a vanishing result for cyclotomic number rings.

\begin{lem}\label{Lem_GaussianReduction}
Let $\zeta_r := e^{2\pi\sqrt{-1}/r}$ be a primitive $r$-th root of
unity for some $r > 0$. Write $r = d\, 2^e$ for some $d > 0$ odd and
$e \geq 0$. Note that
\[
\cO := \Z[\zeta_r] = \Z[\zeta_d, \zeta_{2^e}]
\]
as rings whose involution is complex conjugation. Then, for all $n
\in \Z$, we have
\[
NL_n^h(\cO) = \begin{cases}
NL_n^h(\Z) & \text{if } d=1 \text{ and } e=0,1\\
0 & \text{if } d > 1.
\end{cases}
\]
\end{lem}

\begin{pf}
If $d=1$ and $e=0,1$ then $\cO=\Z$. So we may assume $d>1$.

Write $R := \Z[\zeta_{2^e}]$, and consider $d$ as a divisor of some
odd $N>0$. By the ring decomposition $\rho$ of Proof
\ref{Thm_HyperDecomp} and the isomorphisms of Theorem
\ref{Thm_LocalComplete}, the following induced upper map is an
isomorphism:
\[
\begin{diagram}
\node{\wh{H}^j( R[C_N])} \arrow{e} \arrow{s}
\node{\bigoplus_{d|N} \wh{H}^j( R[\zeta_d])} \arrow{s}\\
\node{\wh{H}^j(R[\textstyle\frac{1}{N}][C_N])} \arrow{e,t}{ \rho_*}
\node{\bigoplus_{d|N} \wh{H}^j( R[\textstyle\frac{1}{N}][\zeta_d]).}
\end{diagram}
\]
But a direct computation shows that the upper map has image in the
$d=1$ factor. That is, the following induced map is an isomorphism:
\[
\incl_*: \wh{H}^j(R[1]) \xra{\quad} \wh{H}^j(R[C_N]).
\]
Then, since $d>1$, the corresponding Tate cohomology groups vanish:
\[
\wh{H}^j(R[\zeta_d]) = 0.
\]
Therefore, since $\cO = R[\zeta_d]$, by Homological Reduction
(\ref{Prop_HomologicalReduction}), we obtain $NL_*^h(\cO)=0$.
\qed\end{pf}

Its analogue in characteristic two is the following lemma.

\begin{lem}\label{Lem_Vanish_Char2Cyclotomic}
Let $P$ be a finite 2-group, and let $d>0$ be odd. Consider the ring
$R=\F_2[P] \xo \Z[\zeta_d]$ with involution. If $d=1$, then for all
$n \in \Z$, the induced map $NL_n(R) \to NL_n(\F_2)$ is an
isomorphism. Otherwise if $d> 1$, then the groups $NL_*(R)$ vanish.
\end{lem}

\begin{pf}
By Nilpotent Reduction (\ref{Prop_NilpotentReduction}), the
following induced map is an isomorphism:
\[
NL_n(R) \xra{\quad} NL_n(\F_2 \xo \Z[\zeta_d]).
\]
Recall, in terms of the $d$-th cyclotomic polynomial $\Phi_d(x) \in
\Z[x]$, that
\[
\F_2 \xo \Z[\zeta_d] = \F_2[x]/(\Phi_d(x)).
\]
Note, by taking formal derivative of $x^d-1$ with $d$ odd, that
$\Phi_d(x)$ is separable over $\F_2$. Then, by the chinese remainder
theorem, the ring $\F_2 \xo \Z[\zeta_d]$ is a finite product of
fields\footnote{It is in fact a product of $\phi(d)/n(d)$ copies of
the finite field $\F_{2^{n(d)}}$, where $\phi(d)$ is the Euler
$\phi$-function and $n(d)>0$ is minimal with respect to the
congruence $2^{n(d)} \equiv 1 \pmod{d}$.} hence is 0-dimensional.

Therefore, by Homological Reducition
(\ref{Prop_HomologicalReduction}), it suffices to show that its Tate
cohomology vanishes. But, as in the previous proof, this follows
from the fact that the induced map
\[
\incl_*: \wh{H}^*( \F_2[1]) \xra{\quad} \wh{H}^*( \F_2[C_N])
\]
is an isomorphism for all odd $N$. \qed\end{pf}

\begin{rem}
It seems appropriate to mention here that the techniques of
Connolly--Ranicki \cite{CR} and of Connolly--Davis \cite{CD} can be
used to generalize their computations of $\UNil_*(\F_2)$.  Namely,
let $\F$ be a perfect field of characteristic two with identity
involution. Here perfect means that the (Frobenius) squaring
endomorphism is surjective. For example, any finite field $\F_{2^e}$
of characteristic two is perfect. Consider the squaring
monomorphisms $\psi$ and $\psi[x]$ defined by
\[
\psi: \F \xra{\quad} \F,\quad  \psi[x]: \F[x] \xra{\quad} \F[x]; \qquad f \longmapsto f^2.
\]

Suppose $n=2k-1$ is odd. Since $\F$ is perfect, by \cite[Proposition
20, Lemma 21, Equation (27)]{CR}, there is an isomorphism
\[
NL_{2k-1}(\F,\id) \xra{\quad} NQ_{2k}(\F,\id) = \Ker\prn{\aug_0: \Ker(\psi[x]-1) \xra{\quad} \Ker(\psi-1)}.
\]
Note that $f \in \F[x]$ and $f^2=f$ imply $f \in \F$. Therefore, we
obtain $NL_{2k-1}(\F,\id)=0$.

Otherwise suppose $n=2k$ is even. The surjectivity of the $\Arf$
invariant below was established earlier by Connolly--Kozniewski
\cite[Proposition 5.7]{CK}.  They asked if $\Arf$ is injective
\cite[Open Question (c)]{CK}.  We answer their question in the
affirmative, as follows.

Consider the $\Z$-module $\F' := \Cok(\psi-1)$. For example, note
$\F_2'=\F_2$. Then, by generalizing the proof of \cite[Lemma 4.6(2),
page 1061]{CD} and using the calculus of \cite[Lemma 4.3]{CD}, it
can be shown that the Arf invariant of symplectic forms over the
function field $\F(x)$ defines an isomorphism
\[
\Arf: NL_{2k}(\F,\id) \xra{\quad} \Cok(\psi[x]-1)/\F' = \bigoplus_{d \text{
odd}} x^d\,\F'.
\]
\end{rem}

The next lemma is a vanishing result for cyclic 2-groups $C$.

\begin{lem}\label{Lem_Vanish_Cyclic}
Let $C$ be a cyclic 2-group, and let $d>1$ be odd. Consider the ring
$R=\Z[\zeta_d]$ whose involution is complex conjugation.  Then the
groups $NL_*(R[C])$ vanish.
\end{lem}

\begin{pf}
We induct on the exponent of $C$. If $e(C)=1$ then
\[
NL_*(R[C]) = NL_*(R) = 0,
\]
by Lemma \ref{Lem_GaussianReduction}(2). Otherwise suppose the lemma
is true for all cyclic 2-groups $C'$ with $e(C') < e(C)$. Then we
may define a ring extension $R'$ of $R$ and a group quotient $C'$ of
$C$ by
\[
R' := R[\zeta_{e(C)}] \quad\text{and}\quad C' := C_{e(C)/2},
\]
and there is a cartesian square (\ref{Thm_Pullback}) of rings with
involution:
\[
\begin{diagram}
\node{R[C]} \arrow{e} \arrow{s} \node{R'} \arrow{s}\\
\node{R[C']} \arrow{e} \node{\F_2[C'] \xo R.}
\end{diagram}
\]
Note, by Lemma \ref{Lem_GaussianReduction}(2) and Lemma
\ref{Lem_Vanish_Char2Cyclotomic}, that
\[
NL_*(R')=0 \quad\text{and}\quad NL_*(\F_2[C'] \xo R)=0.
\]
Therefore, by the Mayer--Vietoris sequence (\ref{Thm_Pullback}), the
map induced by the left column is an isomorphism:
\[
NL_*(R[C]) \xra{\en\iso\en} NL_*(R[C']).
\]
But $NL_*(R[C'])=0$ by inductive hypothesis. This concludes the
argument. \qed\end{pf}

Now we use induction to generalize this vanishing result from cyclic
2-groups $C$ to finite abelian 2-groups $P$.

\begin{lem}\label{Lem_Vanish_Abelian}
Let $R$ be a ring with involution, and let $P$ be a finite abelian
2-group.  If the groups $NL_*(R[C])$ vanish for all cyclic 2-groups
$C$ of exponent $e(C) \leq e(P)$, then the groups $NL_*(R[P])$
vanish.
\end{lem}

\begin{pf}
We induct on the order of $P$.  If $\abs{P}=1$ then
\[
NL_*(R[P]) = NL_*(R[1]) = 0,
\]
by hypothesis.  Otherwise suppose the lemma is true for all $R$ and
$P''$ with $\abs{P''} < \abs{P}$. Since $P$ is a nontrivial abelian
2-group, we can write an internal direct product
\[
P=P' \x C_{e(P)}.
\]
Then we can define a ring extension $R'$ of $R$ and a group quotient
$P_0$ of $P$ by
\[
R' := R[\zeta_{e(P)}] = R[x]/\prn{x^{e(P)/2}+1} \quad\text{and}\quad
P_0 := P' \x C_{e(P)/2}.
\]

Consider the cartesian square (\ref{Thm_Pullback}) of rings with
involution:
\[
\begin{diagram}
\node{R[P]} \arrow{e} \arrow{s} \node{R'[P']} \arrow{s}\\
\node{R[P_0]} \arrow{e} \node{\F_2[P_0] \xo R.}
\end{diagram}
\]
Note, by Nilpotent Reduction (\ref{Prop_NilpotentReduction}) and by
hypothesis using both 2-groups $C$ with $e(C) \leq 2$, that
\[
NL_*(\F_2[P_0] \xo R) \xra{\en\iso\en} NL_*(\F_2 \xo R) = 0.
\]
Also $NL_*(R[P_0]) = 0$, by inductive hypothesis. Then, by the
Mayer--Vietoris sequence (\ref{Thm_Pullback}), the map induced by
the top row is an isomorphism:
\[
NL_*(R[P]) \xra{\en\iso\en} NL_*(R'[P']).
\]
We are done by induction if we can show that $R'$ and $P'$ also
satisfy the hypothesis of the lemma.

Let $C$ be any cyclic 2-group satisfying
\[
1 \leq e(C) \leq e(P') \leq e(P).
\]
We now induct on $e(C)$.  If $e(C)=1$ then
\[
NL_*(R'[C]) = NL_*(R') = 0.
\]
The latter equality follows from the Mayer--Vietoris sequence of the
cartesian square of ring with involution:
\[
\begin{diagram}
\node{R[C_{e(P)}]} \arrow{e} \arrow{s} \node{R'} \arrow{s}\\
\node{R[C_{e(P)/2}]} \arrow{e} \node{\F_2[C_{e(P)/2}] \xo R}
\end{diagram}
\]
and from Nilpotent Reduction, as in the above argument, using the
hypothesis of the lemma.

Otherwise suppose $e(C) > 1$.  Then we may define a quotient group
\[
C' := C_{e(C)/2}
\]
of $C$, and there is a cartesian square of rings with involution:
\[
\begin{diagram}
\node{R'[C]} \arrow{e} \arrow{s} \node{R'[\zeta_{e(C)}]} \arrow{s}\\
\node{R'[C']} \arrow{e} \node{\F_2[C']\xo R'.}
\end{diagram}
\]
We are again done by Nilpotent Reduction and induction on $e(C)$ if
we show that
\[
NL_*(R'[\zeta_{e(C)}])=0.
\]

Consider the primitive root of unity:
\[
\omega := (\zeta_{e(P)})^{e(P)/e(C)} \in R'.
\]
Observe the quotient and factorization
\[
R'[\zeta_{e(C)}] = R[\zeta_{e(P)}][x]/\prn{x^{e(C)/2}+1}
\quad\text{and}\quad x^{e(C)/2} + 1 = \prod_{\text{odd
}d=1}^{e(C)-1} (x - \omega^d).
\]
Then, by the chinese remainder theorem, we obtain an isomorphism of
rings with involution:
\[
R'[\zeta_{e(C)}] \xra{\en\iso\en} \prod_{\text{odd } d=1}^{e(C)-1}
R'.
\]
Hence it induces an isomorphism
\[
NL_*(R'[\zeta_{e(C)}]) \xra{\en\iso\en} \bigoplus_{\text{odd }
d=1}^{e(C)-1} NL_*(R').
\]
But we have already shown that $NL_*(R')=0$. This concludes the
induction on both $e(C)$ and $\abs{P}$. \qed\end{pf}

The last lemma reduces our computation from abelian
2-hyperelementary groups $H$ to abelian 2-groups $P$.

\begin{lem}\label{Lem_Reduction_AbelianHyperelementary}
Consider any abelian 2-hyperelementary group $H=C_N \x P$. Then for
all $n \in\Z$, the following induced map is an isomorphism:
\[
\incl_*: NL_n(\Z[P]) \xra{\quad} NL_n(\Z[H]).
\]
\end{lem}

\begin{pf}
Recall, by Theorem \ref{Thm_HyperDecomp} and additivity of
$L$-groups, that the following induced map is an isomorphism:
\[
NL_n(\Z[H]) \xra{\quad} \bigoplus_{d|N} NL_n(\Z[\zeta_d][P]).
\]
But all the $d\neq 1$ factors vanish by Lemmas
\ref{Lem_Vanish_Cyclic} and \ref{Lem_Vanish_Abelian}. The result now
follows. \qed\end{pf}

We are finally in a position to prove the main theorem.

\begin{pf*}{PROOF of Theorem \ref{Thm_FiniteReduction}.}
Let $H$ be a 2-hyperelementary subgroup of $F$. Since $S$ is normal
abelian, the group $H$ is abelian. Then we can write
\[
H=C_N \x P
\]
for some odd $N$, and $P=H\cap S$ a finite abelian 2-group. So, by
Lemma \ref{Lem_Reduction_AbelianHyperelementary}, the following map
is an isomorphism:
\[
\incl_*: \fN(H\cap S) \xra{\quad} \fN(H).
\]
Therefore, by Hyperelementary Reduction
(\ref{Thm_HyperelementaryReduction}), the following induced map is
an isomorphism:
\[
\incl_*: \fN(S)_{F/S} \xra{\quad} \fN(F).
\]
\qed\end{pf*}

\begin{pf*}{PROOF of Corollary \ref{Cor_OddOrder}.}
This is immediate from the theorem, since $F$ has odd order implies
$S=1$, and since $F/S$ acts trivially by inner automorphisms on $S$.
\qed\end{pf*}
\section{On abelian 2-groups}

Our main theorem (\ref{Thm_FiniteReduction}) reduces the computation
of $\UNil=NL$ for certain finite groups to their maximal abelian
2-subgroup. The following result shows that the $NL$-theory of
abelian 2-groups is determined up to iterated extensions from the
Dedekind domains:
\[
\F_2 \quad\text{and}\quad \Z = \Z[\zeta_2] \quad\text{and}\quad \Z[i] = \Z[\zeta_4]
\quad\text{and}\quad \Z[\zeta_8] \quad\text{and}\quad \Z[\zeta_{16}] \quad\text{and}\quad \ldots.
\]
The $\UNil$-groups of the rings $\F_2$ and $\Z$ with identity
involution and have been calculated \cite{CK, CR, CD, BR}.  The
dyadic cyclotomic number rings $\Z[\zeta_{2^k}]$ for all $k>1$ have
involution given by complex conjugation, and their $\UNil$-groups
shall be calculated in another paper.

\begin{prop}\label{Prop_Abelian2}
Let $P$ be a nontrivial, finite abelian 2-group, and let $n \in \Z$.
Write
\[
P= P' \x C_{e(P)} \quad\text{and}\quad P_0 := P' \x C_{e(P)/2}.
\]

\begin{enumerate}
\item
The Weiss boundary map is an isomorphism:
\[
NQ_{n+1}(\Z[P]) \xra{\en\bdry\en} NL_n(\Z[P]).
\]

\item
There is an exact sequence
\[
\cdots \xra{\quad} NL_{n+1}(\F_2) \xra{\en\bdry\en} NL_n(\Z[P])
\xra{\quad} NL_n(\Z[P_0]) \oplus A \xra{\quad} NL_n(\F_2)
\xra{\en\bdry\en} \cdots
\]
where
\[
A := \begin{cases}NL_n(R_0[P']) & \text{if } P' \neq 1\\
NL_n(\Z[\zeta_{e(P)}]) & \text{if } P'=1 \text{ and } e(P)>1\end{cases} \quad\text{where}\quad R_0 := \Z[\zeta_{e(P)}].
\]

\item
Suppose $R$ is of the form
\[
R = \Z[\zeta_e]
\]
for some $e \geq e(P)$ a power of 2. There is an exact sequence
\[
\cdots \xra{\quad} NL_{n+1}(\F_2) \xra{\en\bdry\en} NL_n(R[P])
\xra{\quad} \bigoplus_2 NL_n(R[P_0]) \xra{\quad} NL_n(\F_2)
\xra{\en\bdry\en} \cdots.
\]
\end{enumerate}
\end{prop}

\begin{pf}
The above sequences are derived from the Mayer--Vietoris exact
sequences (\ref{Thm_Pullback}) of the cartesian squares
\[
\begin{diagram}
\node{\Z[P]} \arrow{e} \arrow{s} \node{R_0[P']} \arrow{s}\\
\node{\Z[P_0]} \arrow{e} \node{\F_2[P_0]}
\end{diagram}
\qquad\text{and}\qquad
\begin{diagram}
\node{R[P]} \arrow{e} \arrow{s} \node{R[P'][\zeta_{e(P)}]} \arrow{s}\\
\node{R[P_0]} \arrow{e} \node{\F_2[P_0].}
\end{diagram}
\]

Part (1) follows by induction on $\abs{P}$ and the five-lemma using
these exact sequences, along with a similar exact sequence for
$\F_2[P]$.  Since $\Z[\zeta]$ and $\F_2$ are Dedekind domains with
involution, the basic cases for the induction are indeed
isomorphisms:
\[
\begin{diagram}
\node{NQ_{n+1}(\Z[\zeta])} \arrow{e,tb}{\bdry}{\iso}
\node{NL_n(\Z[\zeta])} \end{diagram} \qquad\text{and}\qquad
\begin{diagram} \node{NQ_{n+1}(\F_2)} \arrow{e,tb}{\bdry}{\iso}
\node{NL_n(\F_2).}
\end{diagram}
\]

For Parts (2) and (3), recall Nilpotent Reduction
(\ref{Prop_NilpotentReduction}) shows that the following induced map
is an isomorphism:
\[
NL_n(\F_2[P_0]) \xra{\quad} NL_n(\F_2).
\]
Finally, for Part (3), since $\zeta_e \in R$, observe that there
exists an isomorphism of rings with involution:
\[
f: R[\zeta_{e(P)}] \xra{\quad} R[C_{e(P)}]; \qquad \zeta_{e(P)}
\longmapsto (\zeta_e)^{e/e(P)} T,
\]
where $T$ is a generator of the cyclic group $C_{e(P)}$. Therefore
we obtain an induced isomorphism
\[ f_*:
NL_n(R[P'][\zeta_{e(P)}]) \xra{\quad} NL_n(R[P_0]).
\]
\qed\end{pf}

\section{On special 2-groups}

The next step beyond the study of cyclic 2-groups (hence abelian
2-groups, involved in Theorem \ref{Thm_FiniteReduction}) is the
study of special 2-groups. A finite group is \textbf{special} if
every normal abelian subgroup is cyclic.

\begin{prop}[{\cite[2.2.1]{HTW}}]
A finite 2-group $P$ is special if and only if it is either:
\begin{enumerate}
\item[(0)]
for some $e \geq 0$, \textbf{cyclic}
\[
\sC_e := \gens{T \ST T^{2^e}=1}
\]

\item[(1)]
for some $e>3$, \textbf{dihedral}
\[
\sD_e := \gens{T,R \ST T^{2^{e-1}} = 1 = R^2, RTR\inv = T\inv}
\]

\item[(2)]
for some $e>3$, \textbf{semidihedral}
\[
\sSD_e := \gens{T,R \ST T^{2^{e-1}} = 1 = R^2, RTR\inv =
T^{2^{e-2}-1}}
\]

\item[(3)]
for some $e\geq 3$, \textbf{quaternionic}
\[
\sQ_e := \gens{T,R \ST T^{2^{e-1}} = 1, R^2 = T^{2^{e-2}}, RTR\inv =
T\inv}.
\]
\end{enumerate}
\end{prop}

The Mayer--Vietoris exact sequence for cyclic 2-groups $\sC_e$ is
provided in the previous section; the Mayer--Vietoris exact sequence
for the other special 2-groups $P \in \set{\sD_e, \sSD_e, \sQ_e}$ is
provided in the following proposition. The main ingredient is that
$P$ has an index two dihedral quotient $\sD_{e-1}$.

\begin{prop}\label{Prop_Special}
Consider any noncyclic special 2-group $P$, and let $n \in \Z$. Let
$e\geq 3$, and write $\zeta := \zeta_{2^{e-1}}$ a dyadic root of
unity. Denote $\circ_{\pm c}$ as twisting a quadratic extension by
$\pm$ complex conjugation. Then there are the following long exact
sequences.

\begin{multline*}
\cdots \xra{\quad} NL_{n+1}(\F_2) \xra{\en\bdry\en} NL_n(\Z[\sD_e])
\\ \xra{\quad} NL_n(\Z[\sD_{e-1}]) \oplus NL_n(\Z[\zeta] \circ_c C_2)
\xra{\quad} NL_n(\F_2) \xra{\quad} \cdots
\end{multline*}
\begin{multline*}
\cdots \xra{\quad} NL_{n+1}(\F_2) \xra{\en\bdry\en} NL_n(\Z[\sSD_e])
\\ \xra{\quad} NL_n(\Z[\sD_{e-1}]) \oplus NL_n(\Z[\zeta] \circ_{-c}
C_2) \xra{\quad} NL_n(\F_2) \xra{\quad} \cdots
\end{multline*}
\begin{multline*}
\cdots \xra{\quad} NL_{n+1}(\F_2) \xra{\en\bdry\en} NL_n(\Z[\sQ_e])
\\ \xra{\quad} NL_n(\Z[\sD_{e-1}]) \oplus NL_n(\Z[\zeta] \circ_c [i])
\xra{\quad} NL_n(\F_2) \xra{\quad} \cdots
\end{multline*}
\end{prop}

In the next subsections, we shall examine these twisted quadratic
extensions.

\begin{pf}
By Theorem \ref{Thm_Pullback}, we obtain the Mayer--Vietoris exact
sequence
\begin{multline*}
\cdots \xra{\quad} NL_{n+1}(\F_2[G/K]) \xra{\en\bdry\en} NL_n(\Z[G]) \\
\xra{\quad} NL_n(\Z[G/K]) \oplus NL_n(\Z[G]/\Sigma_K) \xra{\quad}
NL_n(\F_2[G/K]) \xra{\quad} \cdots
\end{multline*}
using the following order two subgroup $K$ of the group $G=\sD_e,
\sSD_e, \sQ_e$:
\[
K=\gens{T^{2^{e-2}}}.
\]
Observe that
\begin{gather*}
G/K = \sD_{e-1}\\
\Z[\sD_e]/\Sigma_K = \Z[\zeta] \circ_c C_2\\
\Z[\sSD_e]/\Sigma_K = \Z[\zeta] \circ_{-c} C_2\\
\Z[\sQ_e]/\Sigma_K = \Z[\zeta] \circ_c [i].
\end{gather*}
Finally, Nilpotent Reduction (\ref{Prop_NilpotentReduction}) shows
that
\[
NL_*(\F_2[\sD_{e-1}]) \xra{\quad} NL_*(\F_2)
\]
is an isomorphism. \qed\end{pf}

\subsection{Dihedral and semidihedral 2-groups}

Let $e > 3$ and write $\zeta := \zeta_{2^e}$.  We now show that a
definite chunk of the $\UNil$-groups of the above twisted quadratic
extensions
\[
R := \Z[\zeta] \circ_{\pm c} C_2 = \Z\{\zeta, x\}/(x^2 - 1, x \zeta
x\inv \mp \zeta\inv)
\] consists of the $\UNil$-groups (really $NQ$-groups) of the Dedekind
domains
\[
\cO := \Z[\zeta \pm \zeta\inv].
\]
Here, the involution on $\cO$ is given by complex conjugation.

\begin{prop}\label{Prop_Dihedral_SemiDihedral} For all $n\in \Z$, there is an isomorphism
\[
NL_n(R) \xra{\qquad} NL_n(\cO) \;\oplus\; NL_n\prn{\cO \to R}.
\]
\end{prop}

\begin{pf}
The inclusion $\cO \to R$ of rings with involution induces an exact
sequence of a pair \cite[Proposition 2.2.2]{RanickiExact}:
\[
\cdots \xra{\en\bdry\en} NL_n(\cO) \xra{\en\incl_*\en} NL_n(R)
\xra{\quad} NL_n(\cO \xra{\incl} R) \xra{\en\bdry\en} NL_{n-1}(\cO)
\xra{\quad} \cdots.
\]
Now, W. Pardon \cite[Proof 4.14]{Pardon} constructs an embedding of
rings with involution:
\[ f: R \xra{\quad} M_2(\cO)
\]
whose restriction to the center $\cO$ is the diagonal embedding.
Here, the involution ${}^-$ on the matrix ring $M_2(\cO)$ is defined
by scaling \cite[Defn. 2.5.5]{HTW} the conjugate-transpose
involution ${}^*$:
\[
\ol{B} := A B^* A\inv,
\]
where $A \in SL_2(R)$ is a certain hermitian matrix (i.e. $A =
A^*$). For any quadratic complex $(C,\psi)$ of f.g. projective left
$R$-modules, right-multiplication by $A$ gives an isomorphism from
$\Hom_R(C_i,R)$ with left $R$-module structure given by ${}^*$ to
$\Hom_R(C_i,R)$ with left $R$-module structure given by ${}^-$. Then
we obtain an induced isomorphism
\[
(\id,A)_\#: NL_n(M_2(\cO), {}^*) \xra{\quad} NL_n(M_2(\cO), {}^-).
\]
There is a commutative square
\[
\begin{diagram}
\node{NL_n(\cO)} \arrow{e,t}{\incl_*} \arrow{s,tb}{(\cO\oplus\cO)
\xo_R (-)_\#}{=\diag_*}
\node{NL_n(R)} \arrow{s,b}{f_*}\\
\node{NL_n(M_2(\cO),{}^*)} \arrow{e,t}{(\id,A)_\#}
\node{NL_n(M_2(\cO), {}^-).}
\end{diagram}
\]
The left-hand vertical map is also an isomorphism, by quadratic
Morita equivalence (see \cite[\S 2.4--5]{HTW} for a discussion).
Therefore $\incl_*$ is a split monomorphism, and we obtain the
desired left-split short exact sequence. \qed\end{pf}

\subsection{Quaternionic 2-groups}

Let $e \geq 3$ and write $\zeta := \zeta_{2^e}$. An argument,
similar to the previous subsection, is used to decompose the
$\UNil$-groups of the above twisted quadratic extension
\[
S := \Z[\zeta] \circ_c [i] = \Z\{\zeta, y\} / (y^2 + 1, y \zeta
y\inv - \zeta\inv).
\]
Again consider the Dedekind domain
\[
\cO := \Z[\zeta + \zeta\inv]
\]
with identity involution.

\begin{prop}
For all $n\in \Z$, there is an isomorphism
\[
NL_n(S) \xra{\qquad} NL_n(\cO) \;\oplus\; NL_n\prn{\cO \to S}.
\]
\end{prop}

\begin{pf}
The inclusion $\cO \to S$ of rings with involution induces an exact
sequence of a pair \cite[Proposition 2.2.2]{RanickiExact}:
\[
\cdots \xra{\en\bdry\en} NL_n(\cO) \xra{\en\incl_*\en} NL_n(S)
\xra{\quad} NL_n(\cO \xra{\incl} S) \xra{\en\bdry\en} NL_{n-1}(\cO)
\xra{\quad} \cdots.
\]
By a classical theorem of Weber \cite[Theorem 2.2.4]{HRT}, our ring
$\cO$ of algebraic integers in $F := \Q(\zeta + \zeta\inv)$ is
totally ramified over $2\Z$ by a principal prime ideal $\kp$. That
is, $2 \cO = \kp^r$ where $r := [F : \Q]$.  Performing completions
at $\kp$, W. Pardon \cite[Proof 4.14]{Pardon} constructs an
embedding of rings with involution:
\[ f|: \hat{S}_\kp \subset \hat{\cN}_\kp \xra{\quad} M_2(\hat{\cO}_\kp)
\]
whose restriction to the center $\hat{\cO}_\kp$ is the diagonal
embedding. (There, $\cN \supset S$ is a certain maximal order with
involution in the quaternion algebra $\Q \xo S$ over $F$, so that
$f$ is an isomorphism of $\hat{\cO}_\kp$-algebras with involution.)
Here, the involution ${}^-$ on the matrix ring $M_2(\hat{\cO}_\kp)$
is defined by scaling \cite[Defn. 2.5.5]{HTW} the
conjugate-transpose involution ${}^*$:
\[
\ol{B} := \SmMatrix{0 & 1\\ 1 & 0} B^* \SmMatrix{0 & 1\\ 1 & 0}\inv.
\]
The argument of Proof \ref{Prop_Dihedral_SemiDihedral} shows that
$\incl_*: NL_n(\hat{\cO}_\kp) \to NL_n(\hat{S}_\kp)$ is a split
monomorphism.

Consider the commutative diagram induced by inclusions:
\[
\begin{diagram}
\node{NL_n(\cO)} \arrow{e,t}{\incl_*} \arrow{s} \node{NL_n(S)} \arrow{s}\\
\node{NL_n(\hat{\cO}_{(2)})} \arrow{e} \node{NL_n(\hat{S}_{(2)}).}
\end{diagram}
\]

Observe that the completion of an $\cO$-algebra at $\kp$ equals its
completion at $(2) = 2\cO = \kp^r$.  Then the bottom map is a split
monomorphism.  The vertical maps are isomorphisms by Theorem
\ref{Thm_LocalComplete}(2). Therefore the top map $\incl_*$ is a
split monomorphism, and we obtain the desired left-split short exact
sequence. \qed\end{pf}

\bibliographystyle{elsart-num-sort}
\bibliography{PostDissertation}

\end{document}